\documentclass[a4paper,12pt]{article}
\font\header=cmssdc10 at 20pt

\usepackage[ansinew]{inputenc}   
\usepackage{amsfonts}
\usepackage{indentfirst}
  \usepackage[normalem]{ulem}
  \usepackage{graphicx}
  \usepackage{amssymb}
 \usepackage{draftcopy}
  \usepackage{amsmath}
  \usepackage{rotating}
  \usepackage{color}
  \usepackage{setspace}
  \usepackage{fancybox}
 \usepackage{fancyvrb}
	\usepackage{fancyhdr}
	\usepackage{comment}
	\usepackage[top=0.5cm,left=1cm,right=1cm,bottom=0.5cm]{geometry}
  \usepackage{wrapfig}
    \usepackage{float}
  \usepackage{type1cm}
\usepackage{eso-pic}
\usepackage{color}
 \usepackage{amsthm,amsmath,amsfonts}
\usepackage{fullpage,enumerate}
\usepackage{tikz}
\usepackage[colorlinks,citecolor=blue,urlcolor=blue]{hyperref}
\usepackage{hypernat}
\usepackage{graphicx}

\newcommand{\N}{{\mathbb N}}

\def\1{{\bf 1}}
\def\N{\mathbb N}
\def\Z{\mathbb Z}

\begin{document}

{\header Collective evolution under catastrophes}

\vskip1cm

Rinaldo B. Schinazi

University of Colorado at Colorado Springs

rinaldo.schinazi@uccs.edu

\vskip1cm

{\bf Abstract} We introduce the following discrete time model. Each site of $\N$ represents an ecological niche and is assigned a fitness in $(0,1)$.  All the sites are updated simultaneously at every discrete time. At any given time the environment may be normal with probability $p$ or a catastrophe may occur with probability $1-p$. If the environment is normal the fitness of each site is replaced by the maximum of its current fitness and a random number.  If there is a catastrophe the fitness of each site is replaced by a random number. We compute the joint fitness distribution of any finite number of sites at any fixed time. We also show  convergence of this system to a stationary distribution. This too is computed explicitly.

\bigskip


\vskip1cm

{\header 1 The model}

\vskip1cm

There is strong evidence that the history of the Earth is punctuated by catastrophes: meteor strikes, climate changes, major volcano eruptions and so on. A major catastrophe affects the whole Earth and every ecological niche. The very simple model we introduce follows living species through normal and catastrophic times. We make the following two assumptions. Under normal times the fitness of each species can only increase. Under a catastrophe all the accumulated adaptation is wiped out and the fitness of every species is replaced by a random number. In other words, a catastrophe will provoke a complete renewal of all the ecological niches. Our aim is to study the evolution of such a system.

We now introduce our model. Time is discrete, at each integer $t\geq 0$
each site $n\in\N$ has fitness $\eta_t (n) \in (0,1)$. We may think of each site as an ecological niche. The system of sites evolves in time as follows. Let $p$ be a fixed number in $(0,1)$.  At any time $t\geq 0$ we generate a Bernoulli random variable $B_{t+1}$ with parameter $p$ independent of everything else, and independently a sequence $(U_{t+1}(n):n\in \N)$ of i.i.d. uniform random variables on $(0,1)$. We update the model according to the following rules.
\begin{itemize}
   \item If $B_{t+1}=1$ then  for every $n\in \N$, $\eta_{t+1}(n)  = \max (\eta_t (n), U_{t+1}(n))$.
  \item If $B_{t+1}=0$ then  for every $n\in \N$, $\eta_{t+1}(n)  =  U_{t+1}(n)$.
\end{itemize}

In words, if at time $t$ there is no catastrophe (i.e. $B_{t+1}=1$) then the fitness at time $t+1$ of each site can only go up. On the other hand if at time $t$ there is a catastrophe (i.e. $B_{t+1}=0$) then the fitness at time $t+1$ of each site is reset to a random value.

The joint distribution of $(\eta_t(1),\dots,\eta_t(n))$ for $t\geq 0$ and any natural number $n$ will be shown to depend only on the following function $\phi_t$ defined on $(0,1)$, 

$$
\phi_t(u)
 =u(1-p)\frac{1-(up)^t}{1-up}+u^{t+1}p^t
 %=&\frac{u(1-p)+p^{t+1}u^{t+1}(1-u)}{1-up}
 $$

\medskip

{\bf Theorem 1. }{\sl  Let $(\eta_0(n))_{n\geq 1}$ be an independent sequence of uniform random variables on $(0,1)$. Then, for every $t\geq 0$, $n\geq 1$ and every $(u_1,u_2,\dots,u_n)$ in $(0,1)^n$,
$$ P( \bigcap_{k=1}^n \{\eta_t(k)  \le u_k\})=\phi_t\left(\prod_{k=1}^n u_k\right).$$
}

\medskip

As a consequence of Theorem 1 the cumulative distribution function (c.d.f.) of $\eta_t(n)$ is $\phi_t(u)$ for any $n\geq 1$. Hence, Theorem 1 shows that the c.d.f. of the vector $(\eta_t(1),\dots,\eta_t(n))$ can be expressed using only the c.d.f. at a fixed site. Moreover, Theorem 1 allows explicit computations for the joint distribution of $(\eta_t(1),\dots,\eta_t(n))$.

Let $0<u<1$ and define
\begin{align*}
    \phi(u)=&\lim_{t\to\infty}\phi_t(u)\\
    =&\frac{u(1-p)}{1-up}
\end{align*}

\medskip

{\bf Theorem 2. }{\sl For any initial configuration $\eta_0$, the process $(\eta_t)$ converges in distribution in the following sense. For any $n\geq 1$ and $(u_1,u_2,\dots,u_n)$ in $(0,1)^n$,
$$ \lim_{t\to\infty} P\left( \bigcap_{k=1}^n \{\eta_t(k)  \le u_k\}\right)=\phi\left(\prod_{k=1}^n u_k\right).$$}

\medskip

Note that Theorem 1 holds for a particular initial configuration while Theorem 2 is true for any initial configuration. 
\medskip

{\bf Corollary 1.} {\sl The limiting distribution defined in Theorem 2 is stationary for the process $(\eta_t)$. That is, if $\eta_0$ is distributed according to the limiting distribution so is $\eta_t$ for every $t\geq 0$.}

\medskip

Since a convergence in distribution limit is unique Theorem 2 shows that
the process $(\eta_t)$ has at most one stationary distribution. Corollary 1 shows that in fact the limiting distribution in Theorem 2 is the unique stationary distribution.

Not only does this system converge to a stationary distribution but it does so extremely fast. Figure 1 pictures the histogram of fitness frequencies (fitness on the $x$-axis) after a simulation of the model for 1000 time units. The last catastrophe in this simulation occurred at time 996. At that time we had a flat (uniform) histogram. So it took only 4 steps for the system to "self-organize"! The exponential convergence of this model is also apparent in the expressions for $\phi_t$ and $\phi$. See Figure 2 where we graphed $\phi$ and $\phi_t$ for $t=4$ for the same $p$ we used in the simulation.

It may be interesting to compare this model to the Bak-Sneppen model, see \cite{Bak}. In this discrete time model a finite number of sites are arranged in a circle. At first each site is assigned a random number (i.e. fitness) in $(0,1)$. The system is updated at every discrete time by assigning a new random number to the site with the lowest fitness as well as to its two nearest neighbors. In short, the Bak-Sneppen model evolves only through competition between sites through its "kill the least fit site" (and its unfortunate neighbors) rule.  Such a rule triggers a limiting distribution for which fitnesses below a certain threshold disappear altogether, see also \cite{Ben-Ari1}. In contrast, in our model there is no site competition. Each site fate is entirely driven by its own luck and the environment. Not only that but good times are good for all sites and bad times are bad for all sites. What is remarkable then is how these random events push a flat (uniform) fitness distribution to a distribution with winners and losers. 

With an explicit formula for the limiting distribution we can compute covariances as in the following example. Let $u_1$ and $u_2$ be in $(0,1)$.  Assume that $\eta$ is distributed according to the limiting distribution. For $j=1,2$, let
$$ X_j = {\bf 1}_{\{\eta (j) \le u_j\}}.$$ 
We now compute the covariance of $(X_1,X_2)$.
By Theorem 2,
\begin{align*}
    Cov(X_1,X_2)=&E(X_1X_2)-E(X_1)E(X_2)\\
    =&\phi(u_1u_2)-\phi(u_1)\phi(u_2)\\
    =&(1-p)pu_1u_2\frac{(1-u_1)(1-u_2)}{(1-u_1u_2p)(1-u_1p)(1-u_2p)}
\end{align*}

We see that for all $p$, $u_1$ and $u_2$ in $(0,1)$ this covariance is strictly positive. Therefore, $X_1$ and $X_2$ are positively correlated.

\medskip

Ben-Ari and Schinazi (2022)  have recently considered a similar model where the rule for update is the same as ours under normal times. But under a catastrophe for every $n\in \N$, $\eta_{t+1}(n)  = \min (\eta_t (n),U_{n+1})$. This model turns out to be a lot more difficult to analyze than our model. We will compare the two models in Section 5.

To analyze our model we will follow the general framework of \cite{Ben-Ari}.  But for our model we can take advantage of the renewal aspect (after each catastrophe) of the model. This will allow for a self contained analysis and explicit results at every step.

\vskip1cm

{\header 2 Proof of Theorem 1}

\vskip1cm

{\bf 2.1 Exchangeability}

\bigskip

A sequence of random variables $X_1,X_2\dots$ is said to be exchangeable if for all $n\geq 1$ the vectors $(X_{\sigma(1)},X_{\sigma(2)},\dots,X_{\sigma(n)})$ have the same joint distribution for all permutations $\sigma$  of $\{1,2\dots,n\}$. 
\begin{comment}
A Borel set $B$ is said to be permutable if for all $n\geq 1$ and permutation $\sigma$  of $\{1,2\dots,n\}$,
$$\{(X_1,X_2\dots,X_n)\in B\}=
\{(X_{\sigma(1)},X_{\sigma(2)},\dots,X_{\sigma(n)})\in B\}.$$
The collection of permutable Borel sets turns out to be a $\sigma$-algebra and is called the exchangeable $\sigma$-algebra. 
\end{comment}

Consider an infinite sequence $(Z_n)$ of exchangeable indicators (i.e. a random variable that takes values $0$ and $1$ only) then there exists a random variable $0\leq T\leq 1$ such that
\begin{equation}
P(Z_1=z_1,\dots,Z_n=z_n)=E\left(T^{s_n}(1-T)^{n-s_n}\right),
\end{equation}
where $s_n=z_1+\dots+z_n$, 
see for instance Section 49.3 in Port (1994).

By the symmetry of the dynamics with respect to the sites the stochastic process $(\eta_t)$ is exchangeable in the following sense. If the initial distribution $(\eta_0(n) :n\in\Z_+)$ is exchangeable, for example i.i.d., then for all $t\geq 0$, $(\eta_t(n) :n\in\Z_+)$ is an exchangeable sequence. %We denote the exchangeable $\sigma$-algebra  at time $t$ by ${\cal E}_t$.

We introduce the following sequences of indicators.
Let $u\in (0,1)$ and let 
$$ {\bf I}_t (n,u) = {\bf 1}_{\{\eta_t (n) \le u\}}.$$

Assume that the initial configuration is exchangeable.
Then, for every $t\geq 0$ and $u\in (0,1)$, the sequence $(I_t(n,u):n\in\N)$ is an exchangeable sequence of indicators. 
\begin{comment}
It follows from de Finetti's Theorem that there exists a random variable $\Theta_t (u)$, such that the distribution of  $(I_t(n,u):n\in\N)$ conditioned on ${\cal E}_t$ is i.i.d. with a Bernoulli distribution with parameter  $\Theta_t(u)$. In particular, for every $n\geq 1$, 
$$E(I_t(n,u)|{\cal E}_t)=\Theta_t(u),$$
and therefore 
$$P(\eta_t(n)\leq u)=E(I_t(n,u))=E(\Theta_t(u)).$$
Moreover, we have the following Law of Large Numbers for exchangeable indicators, for fixed $t\geq 0$ and $u\in (0,1)$,
\begin{equation}
\Theta_t (u) = \lim_{N\to\infty} \frac{1}{N} \sum_{n=1}^N {\bf I}_t (n,u), \mbox{ a.s. }
\end{equation}

Furthermore, the exchangeability of the random variables $(\eta_t(n):n\in\N)$ implies that
given ${\cal E}_t$ these random variables are conditionally independent. Hence, for every $t\geq 0$, $n\geq 1$ and every $(u_1,u_2,\dots,u_n)$ in $[0,1]^n$,

\begin{equation}
P( \bigcap_{k=1}^n \{\eta_t(k)  \le u_k\}|{\cal E}_t) = \prod_{k=1}^n P(\eta_t(k)\le u_k|{\cal E}_t) 
=\prod_{k=1}^n \Theta_t (u_k)
\end{equation}
\end{comment}
We apply property (1) to the sequence $(I_t(n,u):n\in\N)$ of exchangeable indicators. We denote the corresponding $T$ in (1) by $\Theta_t(u)$.
\begin{comment}
By (2), for every $t\geq 0$ and $0<u<1$ there exists a random variable $\Theta_t(u)$ such that
\begin{equation}
\Theta_t (u) = \lim_{N\to\infty} \frac{1}{N} \sum_{n=1}^N {\bf I}_t (n,u), \mbox{ a.s. }
\end{equation}
\end{comment}

For fixed $0<u<1$, the process $(\Theta_t(u))_{t\geq 0}$ is updated according to the following rules.
\begin{comment}
This formula shows that the analysis of the process $(\eta_t)$ can be reduced to the analysis of the process $(\Theta_t)$. This is a much simpler Markov process. It is updated according to the following rules.
\end{comment}
\begin{equation}
\Theta_{t+1}(u)  = \begin{cases} u\Theta_t(u)   &\mbox{ if } B_{t+1} = 1 \\ u & \mbox{ if } B_{t+1}=0. \end{cases}
\end{equation}
We now prove this formula. 

If $B_{t+1}=0$ then 
$${\bf I}_{t+1}(n,u)={\bf 1}_{\{U_{t+1}(n)\leq u\}}.$$
The sequence of indicators $\left({\bf 1}_{\{U_{t+1}(n)\leq u\}},n\in\N\right)$ is exchangeable and formula (1) applies. The corresponding $T$ is simply $\Theta_{t+1}(u)=u$. 

On the other hand if $B_{t+1}=1$ then for every $n\in \N$,
$${\bf I}_{t+1}(n,u)={\bf I}_t(n,u){\bf 1}_{\{U_{t+1}(n)\leq u\}}.$$
The sequence of indicators $\left({\bf I}_t(n,u){\bf 1}_{\{U_{t+1}(n)\leq u\}},n\in\N\right)$ is exchangeable. By the independence of $\eta_t$ and $(U_{t+1}(n):n\in \N)$ the corresponding random variable $T$ in (1) is $\Theta_{t+1}(u)=u\Theta_t(u)$. This completes the proof of (2).

We will see below that the distribution of the process $(\eta_t)$ can be computed using the distribution of $(\Theta_t(u))$, a much simpler process.

\bigskip

{\bf 2.2 A renewal process}

\bigskip

Let $T_0=0$ and for $i\geq 1$ let
$$T_i=\min\{s>T_{i-1}:B_s=0\}.$$
That is, $T_i$ is the time of the $i$-th catastrophe. For $t\geq 0$, let 
$$N(t)=\max\{k\geq 0:T_k\leq t\},$$
be the number of catastrophes up to time $t$. It is useful to write $N(t)$ as
$$N(t)=\sum_{i=1}^t B'_i,$$
where $B'_i=1-B_i$. Recall that $(B_i)$ is a sequence of i.i.d. Bernoulli random variables with parameter $p$. Hence,  $(B'_i)$ is a sequence of i.i.d. Bernoulli random variables with parameter $1-p$. This representation of $N(t)$ shows the following two properties for $1\leq s<t$, 

$\bullet$ $N(t)-N(s)$ is independent of $N(s)$.

$\bullet$ $N(t)-N(s)$ has the same distribution as $N(t-s)$.

\medskip

Recall that the cumulative distribution function (c.d.f.) $F$ of a random variable $X$ is defined by $F(x)=P(X\leq x)$.

\medskip

{\bf Proposition 1. }{\sl Let $t>0$, then  $t+1-T_{N(t)}$ has the same distribution as $\min(G,t+1)$ where $G$ has a geometric distribution with parameter $1-p$. That is, the distribution of $G$ is given by $P(G=k)=p^{k-1}(1-p)$ for $k=1,2,\dots$.}

\medskip

{\it Proof}

Since $T_{N(t)}\geq 0$, 
$$t+1-T_{N(t)}\leq t+1.$$
Assume that $1\leq s\leq t+1$.
Then,
\begin{align*}
    P(t+1-T_{N(t)}\geq s)=&\sum_{n=0}^t P(\{N(t)=n\}\bigcap \{t+1-T_n\geq s\}).
\end{align*}

Observe that the event 
$\{t+1-s\geq T_n\})$ is the same as $\{N(t+1-s)\geq n\}$ which can only happen for $t+1-s\geq n$. Hence,
\begin{align*}
 P(t+1-T_{N(t)}\geq s)=&\sum_{n=0}^{t+1-s} P\left(\{N(t)=n\}\bigcap \{N(t+1-s)\geq n\}\right)\\
 =&\sum_{n=0}^{t+1-s}P(\{N(t)-N(t+1-s)=0\}\bigcap \{N(t+1-s)=n\})\\
 =&\sum_{n=0}^{t+1-s}P(N(t)-N(t+1-s)=0) P(N(t+1-s)=n),
\end{align*}
where we used that the random variables $N(t)-N(t+1-s)$ and $N(t+1-s)$ are independent. Since the distributions of $N(t)-N(t+1-s)$ and $N(s-1)$ are the same we get
\begin{align*}
   P(t+1-T_{N(t)}\geq s)=&P(N(s-1)=0)\sum_{n=0}^{t+1-s} P(N(t+1-s)=n)\\
   =&P(N(s-1)=0)\\
   =&p^{s-1},
\end{align*}
for $1\leq s\leq t+1$. By direct computation it is easy to show that
$$P(\min(G,t+1)\geq s)=p^{s-1}.$$
This completes the proof of Proposition 1.

\bigskip

{\bf 2.3 A formula for the underlying Markov chain}

\bigskip

We use the renewal process to get a formula for $\Theta_t(u)$. Let $(\eta_0(n))_{n\geq 1}$ be an independent sequence of uniform random variables on $(0,1)$ then $\Theta_0(u)=u$. For $t\geq 0$,
\begin{equation}
\Theta_t(u)=u^{t+1-T_{N(t)}}.
\end{equation}

\medskip

We now prove (3).
Observe that at times $t=T_n$ for every $n\geq 0$ we have $\Theta_t(u)=u.$
At times $t$ such that $T_n< t<T_{n+1}$ we are strictly in between catastrophes. Hence,
by equation (2)
$$\Theta_t(u)=u\Theta_{t-1}(u).$$
Iterating we get for $T_n\leq t<T_{n+1}$,
$$\Theta_t(u)=u^{t+1-T_n}.$$
Using that $N(t)=n$ if and only if  $T_n\leq t<T_{n+1}$, the preceding equation can be rewritten as
$$\Theta_t(u)=u^{t+1-T_{N(t)}},$$
for all $t\geq 0$. This completes the proof of (3).

We now use (3) to compute the expected value of $\Theta_t(u)$. By Proposition 1, $t+1-T_{N(t)}$ has the same distribution as $\min(G,t+1)$ where $G$ is a geometric random variable with parameter $1-p$. Hence,
\begin{equation*}
P(t+1-T_{N(t)}=k)  = \begin{cases} (1-p)p^{k-1} &\mbox{ if } 1\leq k\leq t \\ p^t & \mbox{ if } k= t+1. \end{cases}
\end{equation*}
Therefore,
\begin{align*}
    E(\Theta_t(u))=&E(u^{t+1-T_{N(t)}})\\
    =&\sum_{k=1}^t u^k (1-p)p^{k-1}+u^{t+1}p^t\\
    =&(1-p)u\frac{1-(up)^t}{1-up}+u^{t+1}p^t\\
    \equiv &\phi_t(u)
\end{align*}

\bigskip

{\bf 2.4 The distribution of the process at a fixed time}

\bigskip

We are now ready to complete the proof of Theorem 1.
Let $(\eta_0(n))$ be an exchangeable sequence.
It follows from de Finetti's Theorem that there exists a $\sigma$-algebra ${\cal E}_t$ such that
the random variables  $(\eta_t(n):n\in\N)$ conditioned on ${\cal E}_t$ are independent,
see Section 57.4 in Port (1994) for instance.
Hence, for every $t\geq 0$, $n\geq 1$ and every $(u_1,u_2,\dots,u_n)$ in $[0,1]^n$,

\begin{equation}
P( \bigcap_{k=1}^n \{\eta_t(k)  \le u_k\}|{\cal E}_t) = 
\prod_{k=1}^n P(\eta_t(k)\le u_k|{\cal E}_t).
\end{equation}
Moreover,
\begin{align*}
 P(\eta_t(k)\le u_k|{\cal E}_t)=&E({\bf 1}_{\eta_t(k)\le u_k}|{\cal E}_t)\\
 =&\Theta_t(u_k).
\end{align*}

By taking expectations across equation (4) we get
$$ P( \bigcap_{k=1}^n \{\eta_t(k)  \le u_k\})=
E\left(\prod_{k=1}^n \Theta_t(u_k)\right).$$

By (3),
\begin{align*}
\prod_{k=1}^n \Theta_t(u_k)=&\left(\prod_{k=1}^n u_k\right)^{t+1-T_{N(t)}}\\
=&\Theta_t\left(\prod_{k=1}^n u_k\right)
\end{align*}

Since $\phi_t(u)=E(\Theta_t(u))$,
by taking expectations on both sides we get
$$ P( \bigcap_{k=1}^n \{\eta_t(k)  \le u_k\})=\phi_t\left(\prod_{k=1}^n u_k\right).$$
This completes the proof of Theorem 1.

\vskip1cm

{\header 3 Proof of Theorem 2}

\vskip1cm

Let $(\tilde\eta_0(n))$ be a sequence of i.i.d. uniform random variables on $(0,1)$. 
Consider now an arbitrary sequence $(\eta_0(n))$ in $(0,1)$, random or deterministic, exchangeable or not. Let $(\eta_t)$ and $(\tilde\eta_t)$ be the processes with initial configurations $\eta_0$ and $\tilde\eta_0$, respectively. We construct $(\eta_t)$ and $(\tilde\eta_t)$ on the same probability space in the following way. At every $t\geq 0$ we use the same Bernoulli $B_{t+1}$ with parameter $p$ and the same sequence $(U_{t+1}(n):n\in \N)$ of uniform random variables to update both processes at time $t+1$. With this construction  we will have for all $t\geq T_1$,
$$\eta_t(n)=\tilde\eta_t(n)\mbox{ for all }n\in \N,$$
where $T_1$ is the time of the first catastrophe (i.e. the first time $t\geq 1$ such that $B_t=0$).
Hence,
\begin{align*}
\left|P( \bigcap_{k=1}^n \{\tilde\eta_t(k)  \le u_k\})-P( \bigcap_{k=1}^n \{\eta_t(k)  \le u_k\})\right|&\leq 2P(\tilde\eta_t(k)\not=\eta_t(k)\mbox{ for some }k)\\
&\leq 2P(T_1>t)\\
&=2p^t.
\end{align*}
By Theorem 1,
$$P( \bigcap_{k=1}^n \{\tilde\eta_t(k)  \le u_k\})=\phi_t\left(\prod_{k=1}^n u_k\right).$$
Since $\phi(u)=\lim_{t\to\infty}\phi_t(u)$ for every $0<u<1$,
$$\lim_{t\to\infty}P( \bigcap_{k=1}^n \{\tilde\eta_t(k)  \le u_k\})=\phi\left(\prod_{k=1}^n u_k\right).$$
Therefore,
$$\lim_{t\to\infty}P( \bigcap_{k=1}^n \{\eta_t(k)  \le u_k\})=\phi\left(\prod_{k=1}^n u_k\right),$$
for any initial configuration $\eta_0$. The proof of Theorem 2 is complete.

\vskip1cm

{\header 4 Proof of Corollary 1}

\vskip1cm

We now prove that the limiting distribution is stationary. Assume that at time $t=0$, $\eta_0$ is distributed according to the limiting distribution. That is, 
for $n\geq 1$ and $(u_1,u_2,\dots,u_n)$ in $[0,1]^n$,
$$ P\left( \bigcap_{k=1}^n \{\eta_0(k)  \le u_k\}\right)=\phi\left(\prod_{k=1}^n u_k\right).$$
By conditioning on the first transition we get,
\begin{align*}
 P\left( \bigcap_{k=1}^n \{\eta_1(k)  \le u_k\}\right)
 =&p P\left( \bigcap_{k=1}^n \{\max(\eta_0(k),U_1(k))  \le u_k\}\right)\\
 +&(1-p)P\left( \bigcap_{k=1}^n \{U_1(k)  \le u_k\}\right)\\
\end{align*}
Using that the random variables $U_1(1),U_1(2),\dots,U_1(n)$ are i.i.d. uniform and independent of the random variables $\eta_0(1),\eta_0(2),\dots,\eta_0(n)$,
\begin{align*}
P\left( \bigcap_{k=1}^n \{\max(\eta_0(k),U_1(k)) \} \le u_k\}\right)=&P\left( \bigcap_{k=1}^n \{\eta_0(k)  \le u_k\}\right)\prod_{k=1}^n u_k\\
=&\phi\left(\prod_{k=1}^n u_k\right)\prod_{k=1}^n u_k.
\end{align*}
Let $u=\prod_{k=1}^n u_k$, we get
\begin{align*}
 P\left( \bigcap_{k=1}^n \{\eta_1(k)  \le u_k\}\right)=&   pu\phi(u)+(1-p)u
\end{align*}
Using now the definition of $\phi$ it is easy to check that
$$pu\phi(u)+(1-p)u=\phi(u).$$
Hence,
$$ P\left( \bigcap_{k=1}^n \{\eta_1(k)  \le u_k\}\right)=\phi\left(\prod_{k=1}^n u_k\right).$$
That is, if $\eta_0$ is distributed according to the limiting distribution so is $\eta_1$. This proves that the limiting distribution is stationary for the process $(\eta_t)$. The proof of Corollary 1 is complete.

\vskip1cm

{\header 5 The underlying Markov chain}

\vskip1cm

In this section we collect results for the Markov chain $(\Theta_t(u))$.

\medskip

$\bullet$ Let $0<u<1$. Assume that $\Theta_0(u)=u$. Then, the Markov chain $(\Theta_t(u))_{t\geq 0}$ converges in distribution, as $t\to\infty$, to $u^G$ where $G$ is a geometric random variable with parameter $1-p$. 

\medskip

We now prove this claim. From Proposition 1, we know that $t+1-T_{N(t)}$ has the same distribution as $G_t=\min(G,t+1)$. As $t$ goes to infinity it is easy to see that $\min(G,t+1)$ converges in distribution to $G$. From formula (3), $(\Theta_t(u))$ has the same distribution as $u^{G_t}$. Hence, $(\Theta_t(u))$ converges in distribution to $u^G$.

\medskip

$\bullet$ Let $0<u<1$. The limiting distribution  of $(\Theta_t(u))_{t\geq 0}$ is stationary.

\medskip

We prove this claim by using generating functions. By conditioning on the first transition,

$$E(s^{\Theta_1(u)})=(1-p)s^u+pE(s^{u\Theta_0(u)}).$$

Assume now that $\Theta_0(u)$ has the same distribution as $u^G$. Then,
\begin{align*}
 E(s^{\Theta_0(u)})&=\sum_{k=1}^\infty (1-p)p^{k-1}s^{u^k}\\
 &=(1-p)s^u+p\sum_{k=2}^\infty (1-p)p^{k-2}s^{u^k}\\
 &=(1-p)s^u+p\sum_{j=1}^\infty (1-p)p^{j-1}s^{u^{j+1}}\\
  &=(1-p)s^u+pE(s^{u\Theta_0(u)})\\
  &=E(s^{\Theta_1(u)})
\end{align*}

This shows that $\Theta_1(u)$ has the same distribution as $\Theta_0(u)$. This completes the proof that $u^G$ is stationary for $(\Theta_t(u))_{t\geq 0}$.

\medskip

$\bullet$ Let $0<u<1$ and $F$ be the cumulative distribution function of $u^G$. Let $0<x<1$. There is a unique $k(x)\in\N$ such that $u^{k(x)}\leq x<u^{k(x)-1}$. Then,
$$F(x)=p^{k(x)-1}.$$

\medskip

Observe that $F$ is a step function with jumps at $u^j$ for all $j\geq 1$. As the steps approach the $x$ axis they are shorter and shorter without ever touching the axis.

\medskip

We now compute $F$.

\begin{align*}
  F(x)=&P(u^G\leq x)\\
  =&\sum_{j\geq 1}P(u^j\leq x)P(G=j)\\
  =&\sum_{j\geq k(x)}P(u^j\leq x)(1-p)p^{j-1}\\
  =&(1-p)\sum_{j\geq k(x)}p^{j-1}\\
  =&p^{k(x)-1}
\end{align*}

\medskip

Closely related to this model is the model introduced in \cite{Ben-Ari}. The dynamics are given by
\begin{itemize}
   \item If $B_{t+1}=1$ then  for every $n\in \N$, $\eta_{t+1}(n)  = \max (\eta_t (n), U_{t+1}(n))$.
  \item If $B_{t+1}=0$ then  for every $n\in \N$, $\eta_{t+1}(n)  =  \min (\eta_t (n), U_{t+1}(n))$.
\end{itemize}
We will call this the $(\max,\min)$ model to differentiate it from our (max,rand) model.
The underlying Markov chain for the $(\max,\min)$ model follows,
\begin{equation}
\Theta_{t+1}(u)  = \begin{cases} \Theta_t(u)  u &\mbox{ if } B_{t+1} = 1 \\ u+(1-u)\Theta_t(u) & \mbox{ if } B_{t+1}=0. \end{cases}
\end{equation}
In \cite{Ben-Ari} it is proved that for the $(\max,\min)$ model,  $\Theta_t(u)$ converges in distribution to 
\begin{equation*}
 \Theta_\infty(u)=\sum_{k=0}^\infty u^{T_k} \left( \frac{1-u}{u}\right)^k,
\end{equation*} 
 where $T_k= G_0 + \dots +G_k$ and
 $G_0,G_1,\dots$ are i.i.d. geometric random variables with parameter $1-p$.  
 Note that the limiting distribution $u^G$ of the (max,rand) model corresponds to the first term (i.e. $k=0$) in the series above. 
 
 Another striking difference between the two models is in the c.d.f. of their limiting distributions. For the $(\max,\min)$ model 
 $\Theta_\infty(u)$ has a c.d.f. which is continuous everywhere but nowhere differentiable. For the (max,rand) model the c.d.f. is differentiable except at the points $u^k$ for all $k\geq 1$. There is, however, a fractal like behavior near 0. 
 
 Fractal like behavior may appear in even simple probability models, see Billingsley (1983) for an interesting example based on the classical ruin problem. Iterated functions systems such as (5) often yield fractals, see Barnsley and Elton (1988) and Strichartz et al. (1995). 
 However, there does not seem to be a clear understanding of why fractals appear. For instance, consider the following iterated function system,
 \begin{equation}
\Theta_{t+1}(u)  = \begin{cases} u\Theta_t(u)   &\mbox{ if } B_{t+1} = 1 \\ 1-u  + u\Theta_t(u) & \mbox{ if } B_{t+1}=0. \end{cases}
\end{equation}
This iterated system has been studied since at least Erdos (1939). There, examples of $u$'s in $(1/2,1)$ are given for which the stationary distribution for the system (6) is continuous but singular with respect to the Lebesgue measure. It is also known that the stationary distribution is absolutely continuous for some values in $(1/2,1)$ and singular for all values in $(0,1/2)$.  As far as we know the question of determining for which $u$'s in $(1/2,1)$ the stationary measure is singular is still open, see also the discussion in \cite[p. 24]{Barnsley}.

 \bibliographystyle{amsplain}

 \begin{figure}[ht] 
\includegraphics[scale=0.5]{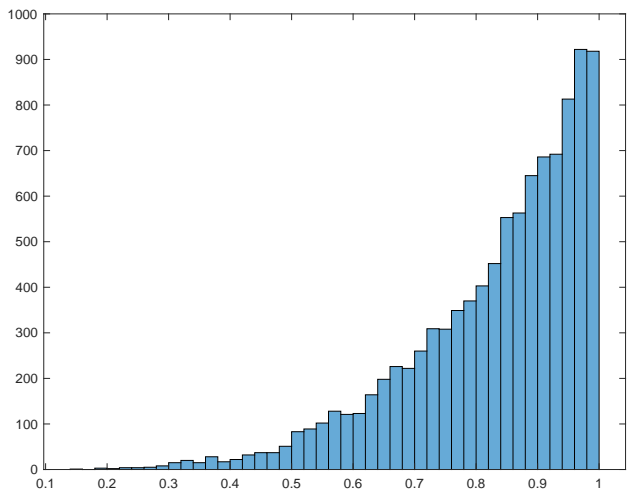}
\caption{This histogram pictures the fitness frequencies (fitness on the $x$-axis) for a simulation which ran for $10^3$ steps for $10^4$ sites and $p=0.9$. The last catastrophe in this simulation occurred at time 996. At that time we had a flat (uniform) histogram.}
\label{fig:1}
\end{figure} 

\begin{figure}[ht] 
\includegraphics[scale=0.5]{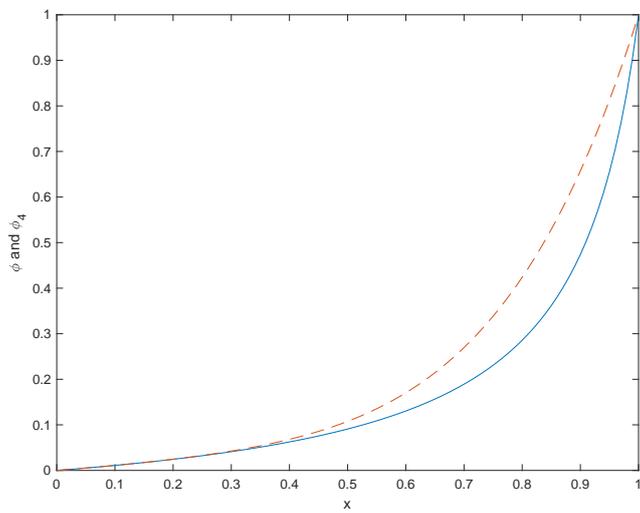}
\caption{ These are the graphs of the cumulative distributions functions $\phi$ and $\phi_4$ at times $t=\infty$ (solid line) and $t=4$ (dashed line) for $p=0.9$. After only 4 time units we see that $\phi_4$ is quite close to $\phi$.}
\label{fig:2}
\end{figure} 
 
\end{document}